\documentclass[a4paper,12pt]{article}

\usepackage[utf8]{inputenc}
\usepackage[ngerman,english]{babel}
\usepackage{amsmath, amsfonts, amssymb, amsthm, mathtools, dsfont, tikz, hyperref}

\usetikzlibrary{arrows.meta}

\allowdisplaybreaks[2]

\makeatletter
\let\@fnsymbol\@alph
\makeatother

\renewcommand{\theenumi}{(\arabic{enumi})}

\def\captionwidth{0.85\linewidth}

\theoremstyle{plain}
\newtheorem{thm}{Theorem}

\theoremstyle{definition}

\newtheorem{rem}[thm]{Remark}

\newcommand{\N}{\mathbb{N}}
\newcommand{\R}{\mathbb{R}}
\newcommand{\C}{\mathbb{C}}
\newcommand{\PP}{\mathbb{P}}
\newcommand{\E}{\mathbb{E}}
\newcommand{\Var}{\operatorname{Var}}
\newcommand{\Cum}{\operatorname{Cum}}
\renewcommand{\Re}{\operatorname{Re}}
\renewcommand{\Im}{\operatorname{Im}}

\newcommand{\gf}{g}					
\newcommand{\cf}{\Psi}				
\newcommand{\quot}{q}				
\newcommand{\nst}{\varrho}			
\newcommand{\nstrad}{r}				
\newcommand{\nstarg}{\varphi}		
\newcommand{\roots}{R}				
\newcommand{\rootsarg}{\Phi}			
\newcommand{\sgn}{\operatorname{sn}}	

\title{Quantification of the Fourth Moment Theorem for Cyclotomic Generating Functions}
\author{%
	Benedikt Redno\ss{}%
	\thanks{Ruhr University Bochum, Germany. Email: benedikt.rednoss@rub.de}
	\and
	Christoph Th\"ale%
	\thanks{Ruhr University Bochum, Germany. Email: christoph.thaele@rub.de}
}
\date{}

\begin{document}
\maketitle

\begin{abstract}
	This paper deals with sequences of random variables $X_n$ only taking values in $\{0,\ldots,n\}$.
	The probability generating functions of such random variables are polynomials of degree $n$.
	Under the assumption that the roots of these polynomials are either all real or all lie on the unit circle in the complex plane, a quantitative normal approximation bound for $X_n$ is established in a unified way.
	In the real rooted case the result is classical and only involves the variances of $X_n$, while in the cyclotomic case the fourth cumulants or moments of $X_n$ appear in addition.
	The proofs are elementary and based on the Stein--Tikhomirov method.
\end{abstract}

\section{Introduction and Results}\label{sec:Intro}


Let $n\in\N$ be a natural number and $X_n$ a random variable that only takes integer values in $\{ 0, \dots , n \}$.
We abbreviate its expectation $\mathbb{E}X_n$ by $\mu_n$,
its variance $\Var(X_n)$ by $\sigma_n^2$
and its fourth cumulant $\Cum_4(X_n) = \E X_n^4-3(\E X_n^2)^2$ by $\kappa_n$. In what follows we shall implicitly assume that $\sigma_n^2 > 0$ to avoid discussion of degenerate cases.
Moreover, by $\gf_n \colon \C \to \C$ we denote the (probability) generating function of $X_n$, that is,
$$
\gf_n(z) := \E[z^{X_n}]= \sum_{k=0}^n \PP(X_n=k) z^k,\qquad z\in\C.
$$
Clearly, $\gf_n$ is a polynomial of degree less than or equal to $n$.
The complex roots of $\gf_n$ will be denoted by $\nst^{(n)}_1, \dots, \nst^{(n)}_n \in \C$, where we do not necessarily assume them to be pairwise distinct.


Starting with the seminal work of Harper \cite{Harper} it has become an active and fruitful direction of research to deduce probabilistic limit theorems for the sequence of random variables $X_n$ as $n\to\infty$ from information about the location of the zeros $\nst^{(n)}_1, \dots, \nst^{(n)}_n$ of $\gf_n$.
Thereby the following two cases have attracted particular interest:
\begin{enumerate}
	\renewcommand{\theenumi}{(\Alph{enumi})}
	\setcounter{enumi}{17}\item \label{case:realrooted}
		$\gf_n$ is \emph{real rooted}, that is, all roots $\nst^{(n)}_1, \dots, \nst^{(n)}_n$ are non-positive real numbers;
	\setcounter{enumi}{2}\item \label{case:cyclotomic}
		$\gf_n$ is \emph{cyclotomic}, meaning that all roots $\nst^{(n)}_1, \dots, \nst^{(n)}_n$ are located on the unit circle in the complex plane.
		In this case, one also says that $g_n$ is \emph{root unitary}.
\end{enumerate}
It is a remarkable observation that a large number of random variables arising in problems having their origin in statistical mechanics or in algebraic and probabilistic combinatorics fall into these two classes.
For example, the number of descents or the number of cycles in a uniform random permutation in finite Coxeter groups (such as the classical groups of type $A$, $B$ or $D$) lead to real rooted generating functions. Also, the matching polynomial associated with an arbitrary graph is real rooted by the celebrated Heilmann--Lieb theorem \cite{HeilmannLieb}. As special cases this yields real rootedness of the monomer-dimer partition function, the $r$-Touchard polynomials or the $r$-Lah number polynomials, for example. We also mention that real rootedness of a polynomial has deep connections with unimodality and log-concavity properties of its coefficient sequence. On the other hand, cyclotomic generating functions arise in the context of the number of inversions in a uniform random permutation in finite Coxeter groups and are prevalent in the theory of polynomials whose coefficients are quotients of so-called $q$-integers. Examples include the $q$-Catalan or the $q$-Narayana numbers, the $q$-analogues of the Hook length formula for standard Young tableaux or the $q$-analogue of the Weyl dimension formula for highest weight modules of semi-simple Lie algebras. Limit theorems for sequences of random variables under assumption \ref{case:realrooted} or \ref{case:cyclotomic} can be found in \cite{Bender} or \cite{HZ15}, respectively.
For more examples, background material and further references we refer the reader to the survey articles \cite{BilleySwanson,Braenden,Pitman}.


Under assumption \ref{case:realrooted} it is a folklore result in probability theory that the sequence of random variables $X_n$ satisfies a central limit theorem if and only if the variances $\sigma_n^2$ satisfy
\begin{align}
	\sigma_n^2\to\infty&\qquad\text{as }n\to\infty,
	\label{iff:realrooted}
\end{align}
see \cite{Bender,Braenden,Pitman}.
On the other hand, in the framework of set-up \ref{case:cyclotomic} it was shown in \cite{HZ15} that a necessary and sufficient condition for $X_n$ to satisfy a central limit theorem is that the fourth cumulants $\kappa_n$ satisfy
\begin{align}
	\frac{\kappa_n}{\sigma_n^4}\to 0&\qquad\text{as }n\to\infty.
	\label{iff:cyclotomic}
\end{align}
The goal of this article is to provide an elementary and unified approach to the following Berry-Esseen type quantitative central limit theorem under assumption \ref{case:realrooted} or \ref{case:cyclotomic}. It is an important observation that the bounds we obtain match the necessary and sufficient conditions presented in \eqref{iff:realrooted} and \eqref{iff:cyclotomic}.
In what follows we let
\begin{align*}
	F(w):=
	\frac{1}{\sqrt{2\pi}}\int_{-\infty}^w e^{-x^2/2}\,\textrm{d}x,\qquad w\in\R,
\end{align*}
be the distribution function of a standard normal random variable.


\begin{thm}\label{thm}
	For $n\in\N$ let $X_n$ be a random variable only taking values in $\{ 0, \dots , n \}$.
	Let $\mu_n=\E X_n$, $\sigma_n^2=\Var(X_n)>0$, $\kappa_n=\Cum_4(X_n)$ and $g_n$ be the generating function of $X_n$. Then, there are absolute constants $c_1>0$ and $c_2>0$ such that the following assertions are valid.
	\begin{enumerate}
		\item\label{thm:i} If $g_n$ is real rooted as defined in \ref{case:realrooted}, it holds that
			\begin{align*}
				\sup_{w\in\R} \Big\vert \PP\Big( \frac{X_n-\mu_n}{\sigma_n} \leq w \Big)-F(w) \Big\vert
				&\leq \frac{c_1}{\sigma_n}.
			\end{align*}
		\item\label{thm:ii} If $g_n$ is cyclotomic as defined in \ref{case:cyclotomic}, one has that
			\begin{align*}
				\sup_{w\in\R} \Big\vert \PP\Big( \frac{X_n-\mu_n}{\sigma_n} \leq w \Big)-F(w) \Big\vert
				&\leq c_2 \frac{\vert\kappa_n\vert^{\frac14}}{\sigma_n}.
			\end{align*}
	\end{enumerate}
\end{thm}


In principle, the result of Theorem~\ref{thm} is known and can be regarded as a special case of the main result of \cite{MichelenSara19Prep}.
In that paper, the roots of the generating function are assumed to be outside a circle of radius $d_n \in (0,1)$ around $1$ or outside a cone pointed at the origin with central angle $\delta_n \in (0,\pi]$, see Figure \ref{fig:delta}.
Under these assumptions, Berry--Essen type bounds of order $\log n/(d_n\sigma_n)$ or $1/(\delta_n \sigma_n)$, respectively, were derived.
Moreover, for random variables whose generating functions are of the form 
\begin{align}
	g_n(z) &= \frac{f(z)}{f(1)}\qquad\text{with}\qquad f(z)=\prod_{j} \frac{1-z^{b_j}}{1-z^{a_j}},\quad a_j,b_j\in\N,
	\label{eq:ass:HST}
\end{align}
a quantitative central limit theorem was obtained in \cite{HeertenEtAl}, among other results, by the method of cumulants.
In this case, the rate of convergence was shown to be of the order $M_n/\sigma_n$, where $M_n:=\max(b_j-a_j)$.
Even though it is known from \cite{HZ15} that in the cyclotomic case \ref{case:cyclotomic} the behaviour of the fourth cumulant regulates whether or not a central limit theorem holds, see \eqref{iff:cyclotomic},
our result is the first to provide a Berry--Esseen type bound including the fourth cumulant. We remark in this context that for all examples treated in \cite{BilleySwanson,HeertenEtAl,HZ15} our result leads to a speed of convergence of the same order and differs only by an absolute constant.
Our computation we present in Section~\ref{sec:proofii} below can in principle be combined with the main result in \cite{MichelenSara19Prep} to derive a result similar to Theorem~\ref{thm}. However, the computations and arguments that are carried out in \cite{MichelenSara19Prep} require substantial theoretical and technical background, though. For this reason, we believe that a self-contained and elementary proof  of  Theorem~\ref{thm} is of independent interest.


\begin{rem}
	Our proof also delivers estimates for the constants $c_1$ and $c_2$ in Theorem~\ref{thm}:
	\begin{align*}
		c_1 \leq 169
		\qquad\text{and}\qquad
		c_2 \leq 447.
	\end{align*}
	This should be compared to the constants $2^{3261}$ and $2^{3257}$ from \cite{MichelenSara19Prep} for the two already mentioned more general cases,
	and to the constant $324\sqrt{2} \cdot (\pi \sqrt{7/24})^{-1} \leq 271$ from \cite{HeertenEtAl} for the special case satisfying \eqref{eq:ass:HST}.
\end{rem}


As we already pointed out above, our contribution is an elementary and unified route to Theorem~\ref{thm},
which is based on the so-called Stein--Tikhomirov method.
The classical Stein's method for normal approximation is based on the observation
that a real-valued random variable $Z$ is standard normally distributed
if and only if
\begin{equation}\label{eq:SteinCharcaterization}
	\E[ f^\prime(Z) - Zf(Z) ]=0
\end{equation}
for all sufficiently regular functions $f \colon \R\to\R$.
If $Z$ is not standard normally distributed, the expectation $\E[ f^\prime(Z) - Zf(Z) ]$ can be used to quantify how `close' the distribution of $Z$ is to the standard normal distribution.
Eventually, this leads to bounds for common probability distances,
such as the Kolmogorov distance or the Wasserstein distance.
For a detailed introduction to Stein's method, see \cite{CGS11}.
Tikhomirov \cite{Ti80} combined this core idea of Stein's method with the more traditional theory of characteristic functions.
Choosing for fixed $t\in\R$ the function $f$ in \eqref{eq:SteinCharcaterization} as $f(x) = e^{itx}$, $x \in \R$, leads to the ordinary differential equation
\begin{align}
	t \cf_Z(t) + \cf_Z^\prime(t) = 0 \label{eq:STM}
\end{align}
for all $t \in \R$,
where $\cf_Z(t) := \E[e^{itZ}]$ denotes the characteristic function of $Z$.
The differential equation in \eqref{eq:STM} indeed characterizes the characteristic function of the standard normal distribution, which is $e^{-\frac12 t^2}$, $t \in \R$.
This motivates the idea of evaluating the left-hand side of \eqref{eq:STM} even if $Z$ is not exactly standard normally distributed,
in order to derive a non-uniform bound for the difference $\vert \cf_Z(t) - e^{-\frac12 t^2} \vert$ of characteristic functions.
This bound can then be combined with Esseen's classical smoothing lemma, which states that
\begin{align*}
	\sup_{w\in\R} \big\vert \PP(Z \leq w) - F(w) \big\vert
				&\leq \frac1\pi \int_{-T}^T \bigg\vert \frac{\cf_Z(t) - e^{-\frac12 t^2}}{t} \bigg\vert\,\textrm{d}t + \frac{12\sqrt{2}}{\sqrt{\pi^3}T}
\end{align*}
for all $T>0$, see \cite[Lemma 2 in Chapter XVI.3]{Feller2}. Many important aspects of this strategy are carried out rather implicitly in the original paper \cite{Ti80},
but a lemma summarizing and partly generalizing the essential ingredients of this approach is the content of \cite[Lemma 2.3]{Ro22}. A version of it can be rephrased as follows:
Suppose there are continuous functions $a,b \colon \R \to \C$ and non-negative constants $A_0, A_1, B_0, B_1, B_2,C \in \R$
satisfying $A_0 < \frac12$,
$\cf_Z^\prime(t) = -t (1 + a(t)) \cf_Z(t) + b(t)$,
$\vert a(t) \vert \leq A_0 + A_1 \vert t \vert$, and
$\vert b(t) \vert \leq B_0 + B_1 \vert t \vert + B_2 t^2$
for all $t\in(-C,C)$.
Then,
\begin{align}
	\nonumber\sup_{w\in\R} \big\vert \PP(Z \leq w) - F(w) \big\vert
	\leq{}&
	\frac2\pi A_0
	+ \frac{4}{3\sqrt{\pi}} A_1
	+ \sqrt{\pi} B_0 \\
	\nonumber&+ \frac2\pi \big( 1 + 2 \max \{ 0, -\log 2s \} \big) B_1 \\
	&+ \frac4{\pi s} B_2
	+ \frac{24s}{\pi \sqrt{2\pi}} \label{eq:SteinTikBound}
\end{align}
for all $s \geq \max \big\{ \frac{2A_1}{1-2A_0} , \frac{1}{C} \bigr\}$. In fact, this follows from Lemma~2.3 in \cite{Ro22} by realizing that the assumptions on $a(t)$ and $b(t)$ in said lemma actually are only required to hold in an interval $(-C,C)$ around zero and not necessarily for the whole real line.

The Berry-Esseen type bound  \eqref{eq:SteinTikBound} is the starting point of our proof of Theorem \ref{thm}, which in turn is the content of the remaining parts of this paper. More precisely, in Section \ref{sec:PrelimComp} we collect some preliminary computations and introduce some further notation. The proof of Theorem \ref{thm} is the content of the final Section \ref{sec:Proof}. 

\section{Preliminary Computations}\label{sec:PrelimComp}


To simplify our presentation, in this section and the next one we fix $n\in\N$ and suppress the index $n$ in what follows.
In particular, we write
$X$ for a random variable $X_n$ taking values in $\{0,\ldots,n\}$,
$\mu$ for the expectation,
$\sigma^2$ for the variance,
$\kappa$ for the fourth cumulant,
$\gf$ for the generating function of $X$
and
$\nst_1, \dots, \nst_n$ for the (not necessarily distinct) complex roots of $g$. We also put $\roots := \{ \nst_1, \dots , \nst_n \}$.


We start by noting that for all $z \in \C$,
\begin{align}
	\gf(z)
		&= \E[z^X]
		= \sum_{k=0}^n \PP(X=k) z^k\nonumber\\
		&= \PP(X=n) \prod_{k=0}^n ( z - \nst_k )
		, \nonumber
		\intertext{and}
	\gf^\prime(z)
		&= \PP(X=n) \sum_{k=0}^n \prod_{\substack{ \ell=0 \\ \ell \neq k}}^n ( z - \nst_\ell )
		.\label{gf-prime}
\end{align}
Since all coefficients of $\gf$ are real-valued and positive or zero, there cannot be any root on the positive real axis $\{z\in\C:\Im(z)=0,\Re(z)>0\}\subset\C$.
Further, if there is a root with non-zero imaginary part, then its complex conjugate is a root as well and it also has the same multiplicity.
Excluding the case of $z$ being a root of $\gf$, we can rewrite \eqref{gf-prime} as follows:
\begin{align}
	\gf^\prime(z)
		&= \gf(z)  \sum_{k=0}^n \frac1{z - \nst_k}
	\label{gf-prime-relation-to-gf}
\end{align}
for all $z \in \C \setminus \roots$.

\begin{figure}
	\centering
	\begin{tikzpicture}[scale=2]
		\def\deltavalue{40}
		\def\deltaradius{2}
		\path ({180+\deltavalue}:\deltaradius) coordinate (LU) -- (\deltavalue:\deltaradius) coordinate (RO);
		\path  (0,0) coordinate (CC)
			-- (LU |- RO) coordinate (LO)
			-- (RO |- LU) coordinate (RU)
			-- (LU |- CC) coordinate (LC)
			-- (RO |- CC) coordinate (RC)
			-- (CC |- RO) coordinate (CO)
			-- (CC |- LU) coordinate (CU)
			;
		\fill[color=green!8!white,domain=0:1.9,samples=100] (CC) -- (\deltavalue:\deltaradius) -- (LO) -- (LU) -- (-\deltavalue:\deltaradius) -- cycle;
		\fill[color=red!8!white,domain=0:1.9,samples=100] (CC) -- (\deltavalue:\deltaradius) -- (-\deltavalue:\deltaradius) -- cycle;
		\draw[Stealth-] (RC) ++ (0.2,0) -- (LC);
		\draw[Stealth-] (CO) ++ (0,0.2) -- (CU);
		\begin{scope}
			\clip (LU) rectangle (RO);
			\draw[color=orange,very thick] (RU) -- (CC) -- (RO);
			\draw (0.5,0) arc (0:\deltavalue:0.5);
			\path ({\deltavalue/2}:0.33) node[font=\scriptsize] {$\delta$};
			\draw[color=blue,very thick] (CC) circle(1);
			\draw[color=black,dashed] (0,0) -- (200:1) node[pos=0.5,sloped,below,font=\scriptsize,inner sep=2] {$\vert z\vert=1$};
			\foreach \a in {
					(  \deltavalue:1.1),
					( -\deltavalue:1.1),
					(  50:1.0),
					( -50:1.0),
					(  53:1.4),
					( -53:1.4),
					(  65:1.2),
					( -65:1.2),
					(  80:0.8),
					( -80:0.8),
					( 105:1.1),
					(-105:1.1),
					( 120:1.4),
					(-120:1.4),
					( 125:0.7),
					(-125:0.7),
					( 135:1.6),
					(-135:1.6),
					( 150:0.8),
					(-150:0.8),
					( 160:1.3),
					(-160:1.3),
					( 180:1.4)
				}
			\path[fill=green!75!black,radius=0.03] \a circle;
		\end{scope}
	\end{tikzpicture}
	\begin{minipage}{\captionwidth}
		\caption{Roots of the generating function $\gf$.
			The angle parameter $\delta$ is chosen such that no root (green dots) of $\gf$ is in the interior of the cone pointed at the origin and having central angle $\delta$ (red).
			Evaluation of $\gf$ around the unit circle (blue) yields the characteristic function $\cf$.
			Note that we currently do not necessarily assume $\gf$ to be real rooted or cyclotomic as in~\ref{case:realrooted} or~\ref{case:cyclotomic}, respectively.}
		\label{fig:delta}%
	\end{minipage}%
\end{figure}
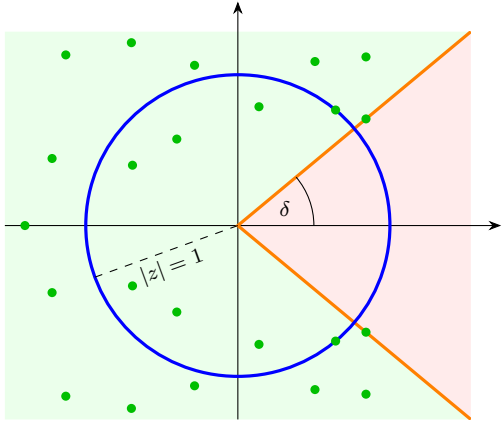

For every $k\in\{ 1, \dots , n \}$,
let $(\nstrad_k,\nstarg_k) \in [0, \infty) \times (-\pi, \pi]$
be the polar coordinates of $\nst_k$,
i.e.\ $\nst_k = \nstrad_k  e^{i\nstarg_k}$.
If $0 \in \C$ is a root, we assign to it the polar coordinates $(0,\pi)$.
We define the angle set $\rootsarg := \{ \nstarg_1, \dots , \nstarg_n \}$
and denote $\delta := \min_{\nstarg \in \rootsarg} \vert\nstarg\vert$.
Since $0 \not\in \rootsarg$, we have that $\delta>0$, see Figure~\ref{fig:delta} for an illustration.


Let $\cf \colon \R \to \C$ be the characteristic function of $X$ given by $\cf(t) := \E[e^{itX}]$, $t\in\R$.
In particular, $\cf$ arises from the generating function $\gf$ via the change of variables $t\mapsto e^{it}$,
so that
\begin{align}
	\cf(t)
		&= g(e^{it})
		\nonumber
		\intertext{and}
	\cf^\prime(t)
		&= i  e^{it} \, g^\prime(e^{it})
		\label{cf-prime}
\end{align}
for all $t \in \R$.
Similar to \eqref{gf-prime-relation-to-gf},
relation \eqref{cf-prime} can be rewritten as
\begin{align*}
	\cf^\prime(t)
		&= i \sum_{k=0}^n \frac{e^{it}}{e^{it} - \nst_k}\, g(e^{it})
\end{align*}
for all $t \in \R \setminus \rootsarg$.
Therefore, defining the function $\quot \colon \R \setminus \rootsarg \to \C$ by
\begin{align}\label{eq:defQ}
	\quot(t) := i \sum_{k=0}^n \frac{e^{it}}{e^{it} - \nst_k},
\end{align}
we arrive at
\begin{align*}
	\cf^\prime(t)
		&= \quot(t) \, g(e^{it}) \\
		&= \quot(t) \, \cf(t) ,
		\intertext{and}
	\cf^{\prime\prime}(t)
		&= \quot^\prime(t) \, \cf(t) + \quot(t) \, \cf^\prime(t)\\
		&= \bigl( \quot^\prime(t) + (\quot(t))^2 \bigr) \, \cf(t)
\end{align*}
for all $t \in \R \setminus \rootsarg$.


In preparation for later calculations,
we compute the first three derivatives of $\quot$.
For $t \in \R \setminus \rootsarg$ we have
\begin{align}
	\quot^\prime(t)
		&= i\sum_{k=0}^n \frac{ie^{it}(e^{it} - \nst_k) - e^{it} ie^{it}}{(e^{it} - \nst_k)^2}\nonumber\\
		 &=  \sum_{k=0}^n \frac{ \nst_k e^{it} }{(e^{it} - \nst_k)^2}
		\nonumber \\
		&=  \sum_{k=0}^n \biggl(
				\frac{ -1 }{1 - \nst_k e^{-it}}
				+ \Bigl( \frac{ -1 }{1 - \nst_k e^{-it}} \Bigr)^2
			\biggr)
		\nonumber
\intertext{and}
	\quot^{\prime\prime}(t)
		&= \sum_{k=0}^n \biggl(
				\frac{ (-\nst_k)(-ie^{-it}) }{(1 - \nst_k e^{-it})^2}
				+ 2
				\cdot \frac{ - 1 }{1 - \nst_k e^{-it}}
				\cdot \frac{ (-\nst_k)(-ie^{-it}) }{(1 - \nst_k e^{-it})^2}
			\biggr)
		\nonumber \\
		&= i \sum_{k=0}^n
				\frac{ \nst_k e^{-it} }{(1 - \nst_k e^{-it})^2}
				\biggl( 1 - \frac{ 2 }{1 - \nst_k e^{-it}} \biggr)
		,\label{quot-primeprime}
\intertext{whereas the third derivative satisfies}
	\quot^{\prime\prime\prime}(t)
		&= i \sum_{k=0}^n \biggl(
				\frac{
					-i \nst_k e^{-it} (1 - \nst_k e^{-it})^2
					- \nst_k e^{-it} 2 (1 - \nst_k e^{-it}) i\nst_k e^{-it}
				}{(1 - \nst_k e^{-it})^4}
		\nonumber\\
		&\hspace{2cm}
				-2\frac{
					-i \nst_k e^{-it} (1 - \nst_k e^{-it})^3
					- \nst_k e^{-it} 3 (1 - \nst_k e^{-it})^2 i\nst_k e^{-it}
				}{(1 - \nst_k e^{-it})^6}
			\biggr)
		\nonumber \\
		&= \sum_{k=0}^n
				\frac{ \nst_k e^{-it} }{(1 - \nst_k e^{-it})^3}
				\biggl( 1 + \nst_k e^{-it} - 2 \frac{ 1 + 2 \nst_k e^{-it} }{1 - \nst_k e^{-it}} \biggr)
		\nonumber \\
		&= \sum_{k=0}^n
				\frac{ -\nst_k e^{-it} }{(1 - \nst_k e^{-it})^4}
				\Bigl( (1 + \nst_k e^{-it})^2 + 2 \nst_k e^{-it} \Bigr)
		.\nonumber
\end{align}

\begin{figure}
	\centering
	\begin{tikzpicture}[scale=1.5]
		\path (-3.5,-1.8) coordinate (LU) -- (4.0,1.8) coordinate (RO);
		\path  (0,0) coordinate (CC)
			-- (LU |- RO) coordinate (LO)
			-- (RO |- LU) coordinate (RU)
			-- (LU |- CC) coordinate (LC)
			-- (RO |- CC) coordinate (RC)
			-- (CC |- RO) coordinate (CO)
			-- (CC |- LU) coordinate (CU)
			;
		\begin{scope}
			\clip (LU) rectangle (RO);
			\fill[color=green!8!white,domain=0:1.9,samples=100]
				plot ({\x},{ sqrt(\x*(\x-1)^2/(2-\x)}) -- (LO) -- (LC) -- cycle;
			\fill[color=green!8!white,domain=0:1.9,samples=100]
				plot ({\x},{-sqrt(\x*(\x-1)^2/(2-\x)}) -- (LU) -- (LC) -- cycle;
			\fill[color=red!8!white,domain=0:1.9,samples=100]
				plot ({\x},{ sqrt(\x*(\x-1)^2/(2-\x)}) -- (RO) -- (RC) -- cycle;
			\fill[color=red!8!white,domain=0:1.9,samples=100]
				plot ({\x},{-sqrt(\x*(\x-1)^2/(2-\x)}) -- (RU) -- (RC) -- cycle;
		\end{scope}
		\draw[Stealth-,overlay] (RC) ++ (0.2,0) node[below] {$x$} -- (LC);
		\draw[Stealth-] (CO) ++ (0,0.2) node[right] {$y$} -- (CU);
		\begin{scope}
			\clip (LU) rectangle (RO);
			\draw[color=orange,very thick,domain=0:1.9,samples=100]
				plot ({\x},{sqrt(\x*(\x-1)^2/(2-\x)})
				plot ({\x},{-sqrt(\x*(\x-1)^2/(2-\x)});
			\path (1.7,0 |- CU) node[above left,font=\footnotesize,inner sep=2,color=orange] {$y^2 = \frac{x (x-1)^2}{2-x}$};
			\draw[color=blue,very thick]
				(LC) -- (CC) circle(1) -- (RC);
			\path (-1,0) node[above left,font=\scriptsize,inner sep=2,color=blue] {$y  (\vert\nst\vert-1)=0$};
			\draw[color=black,dashed]
				(2,0 |- CU) -- (2,0 |- CO) node[pos=0.25,sloped,below,font=\scriptsize,inner sep=2] {$x=2$};
			\draw[color=black,dashed]
				(0,0) -- (200:1) node[pos=0.5,sloped,below,font=\scriptsize,inner sep=2] {$\vert\nst\vert=1$};
		\end{scope}
		\path (LO) -- (CC) node[pos=0.5,above,color=green!66!black] {$\Re\bigr( \frac{-\nst}{(1-\nst)^2} \bigr) \geq 0$};
		\path (RO) -- (2,0) node[pos=0.5,above,color=red] {$\Re\bigr( \frac{-\nst}{(1-\nst)^2} \bigr) < 0$};
		\path (225:1) node[below left,color=blue] {$\Im\bigr( \frac{-\nst}{(1-\nst)^2} \bigr) = 0$};
	\end{tikzpicture}
	\begin{minipage}{\captionwidth}
		\caption{Behaviour of $\frac{-\nst}{(1-\nst)^2}$;
			$x$ denotes $\Re(\nst)$, $y$ denotes $\Im(\nst)$.
			According to \eqref{sigma-re},
			the real part of $\frac{-\nst}{(1-\nst)^2}$ is positive or zero for all $\nst$ in the green area.
			According to \eqref{sigma-im},
			the imaginary part of $\frac{-\nst}{(1-\nst)^2}$ is zero for all $\nst$ on the blue line and the blue circle.
			The orange curve is the
			\emph{right strophoid}
			of the curve given by $x=2$
			with respect to the pole $0+0i$
			and the fixed point $1+0i$,
			see page~96 in \cite{Lo61}.}
		\label{fig:sophoide}%
	\end{minipage}
\end{figure}
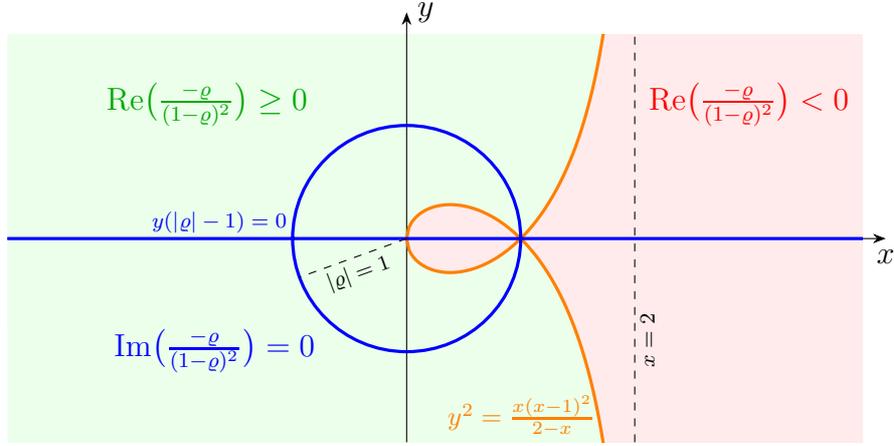

By definition,
the function $\quot$ is closely related to the cumulants of $X$.
In particular, for the first, the second and the fourth cumulant we have that
\begin{align}
	\mu
		&= \E[X]
		= -i \cf^\prime(0)
		= -i \quot(0) \, \cf(0) \nonumber\\
		&= -i \quot(0)
		= \sum_{k=0}^n \frac{1}{1-\nst_k}
		, \label{expectation}\\
	\sigma^2
		&= \Var(X)
		= - \cf^{\prime\prime}(0) - \mu^2
		= - \bigl( (\quot(0))^2 + \quot^\prime(0) \bigr) + (\quot(0))^2 \nonumber\\
		&= - \quot^\prime(0)
		= \sum_{k=0}^n \frac{-\nst_k}{(1-\nst_k)^2}
		\label{variance}
		\intertext{and}
	\kappa
		&= \Cum_4(X) = \quot^{\prime\prime\prime}(0) \nonumber\\
		&= \sum_{k=0}^n
				\frac{ -\nst_k }{(1 - \nst_k )^4}
				\bigl( (1 + \nst_k )^2 + 2 \nst_k \bigr)
		. \label{fourth-cumulant}
\end{align}
Regarding the variance in \eqref{variance},
it is interesting to note that for any root $\nst\in\roots$,
\begin{align*}
	\frac{-\nst}{(1-\nst)^2}
		&= \frac{-\nst \cdot (1-\bar\nst)^2}{(1-\nst)^2(1-\bar\nst)^2}
		= \frac{-\nst + 2 \vert\nst\vert^2 - \bar\nst\vert\nst\vert^2}{\vert 1-\nst \vert^4}
	.
\end{align*}
Straightforward calculations lead now to the crucial observation that
\begin{align}
	\Re\Bigr( \frac{-\nst}{(1-\nst)^2} \Bigr) &\geq 0
	& \iff &&
	&\Im(\nst)^2 \geq \frac{\Re(\nst) (\Re(\nst)-1)^2}{2-\Re(\nst)},
	\label{sigma-re}\\
	\Im\Bigr( \frac{-\nst}{(1-\nst)^2} \Bigr) &= 0
	& \iff &&
	&\Im(\nst)=0 \text{ or } \vert\nst\vert=1
	.\label{sigma-im}
\end{align}
Therefore, the expression $\frac{-\nst}{(1-\nst)^2}$ is real and positive or zero if and only if $\nst$ is real and negative or zero or $\nst$ has modulus $1$,
see Figure~\ref{fig:sophoide}.
These are precisely the cases \ref{case:realrooted} and \ref{case:cyclotomic} we consider in Theorem~\ref{thm}.

\section{Proof of Theorem~\ref{thm}}\label{sec:Proof}

\subsection{Common Ground}


We use the same notation as in the previous section and let $X$ be a random variable taking values in $\{0,\ldots,n\}$ for some $n\in\N$, $\mu$ and $\sigma^2>0$ stand for its expectation and variance, respectively.
Let $\widehat X$ be the standardization of $X$, that is,
$$
\widehat X = \frac{1}{\sigma} (X - \mu).
$$
We write $\widehat \cf$ for the characteristic function of $\widehat X$.
Then, 
\begin{align*}
	\widehat\cf(t)
		&= \cf \Bigl( \frac{t}{\sigma} \Bigr) \, e^{-i\frac{\mu}{\sigma}t}
		\intertext{and}
	\widehat\cf^\prime(t)
		&= \frac1\sigma
			\, \cf^\prime\Bigl( \frac{t}{\sigma} \Bigr)
			\, e^{-i\frac{\mu}{\sigma}t}
			- i \frac\mu\sigma
			\,\cf\Bigl( \frac{t}{\sigma} \Bigr)
			\, e^{-i\frac{\mu}{\sigma}t}
\end{align*}
for all $t \in \R$.
Using the function $\quot$ defined by \eqref{eq:defQ} we find
\begin{align*}
	\widehat\cf^\prime(t)
		&= \frac1\sigma
			\, \Bigl( \quot\Bigl( \frac{t}{\sigma} \Bigr) - i\mu \Bigr)
			\, \widehat\cf(t)
\end{align*}
for $t \in \R\setminus\sigma\rootsarg$ with $\sigma\rootsarg=\{\sigma\varphi:\varphi\in\Phi\}$.


Following the Stein--Tikhomirov method as explained in Section \ref{sec:Intro}, we are interested for $t \in \R \setminus \sigma\rootsarg$ in the differential expression
\begin{align*}
	\widehat\cf^\prime(t) + t \widehat\cf(t)
	&= \biggl( \frac{ \quot\bigl( \frac{t}{\sigma} \bigr) - i\mu }\sigma + t \biggr) \,\widehat\cf(t)
	= \alpha(t) \, \widehat\cf(t)
\end{align*}
with the continuous function
\begin{align*}
	\alpha(t) :={}& \frac{ \quot\bigl( \frac{t}{\sigma} \bigr) - i\mu }\sigma + t.
\end{align*}
We aim to show that under the assumptions \ref{case:realrooted} or \ref{case:cyclotomic} on the location of the roots of $\gf$ there are constants $A > 0$ and $C > 0$ such that
\begin{align}
	\vert \alpha(t) \vert &\leq A \, t^2 \qquad \text{for all}\qquad t \in (-C,C).
	\label{assumption}
\end{align}
We can then apply the Stein--Tikhomirov bound \eqref{eq:SteinTikBound} with the functions $a(t)$ and $b(t)$ there given by $-ta(t)=\alpha(t)$ and $b(t)\equiv 0$.
This leads to $A_0=B_0=B_1=B_2=0$ and $A_1=A$ and with the choice $s=\max\bigl\{ 2A, \frac1C \bigr\}$, \eqref{eq:SteinTikBound} then leads to
\begin{align}
	\sup_{w \in\R} \vert \PP(\widehat X \leq w) - F(w)\vert
	&\leq
	\frac{4}{3\sqrt\pi} \, A
	+ \frac{24}{\pi\sqrt{2\pi}} \, \max\Bigl\{ 2A, \frac1C \Bigr\}.
	\label{Stein-Tikho}
\end{align}


In order to verify \eqref{assumption}
we use \eqref{expectation} to write $i\mu = \quot(0)$ and
we use \eqref{variance} to see that
$t = -\frac{1}{\sigma} \cdot \quot^\prime(0) \cdot \frac{t}{\sigma}$.
Then,
\begin{align*}
	\alpha(t)
	&= \frac{
			\quot\bigl( \frac{t}{\sigma} \bigr)
			- \quot(0)
			- \quot^\prime(0) \cdot \frac{t}{\sigma}
		}{\sigma}
\end{align*}
for $t \in \bigl( -\frac{\sigma\delta}{2} , \frac{\sigma\delta}{2} \bigr) \subset \R \setminus \sigma\rootsarg$.
This equals the remainder term of a Taylor polynomial, which implies that $\alpha(t)$ satisfies the estimate
\begin{align}
	\vert \alpha(t) \vert
	&\leq \frac{1}{2\sigma}
		\sup_{s \in ( -\frac{\sigma\delta}{2} , \frac{\sigma\delta}{2} )}
		\Bigl\vert \quot^{\prime\prime} \Bigl( \frac{s}{\sigma} \Bigr) \Bigr\vert
		\cdot
		\frac{t^2}{\sigma^2}
	= \frac{t^2}{2\sigma^3}
		 \sup_{s \in ( -\frac{\delta}{2} , \frac{\delta}{2} )}
		\vert \quot^{\prime\prime}(s) \vert
	\label{alpha-estimate}
\end{align}
for all $t \in ( -\frac{\sigma\delta}{2} , \frac{\sigma\delta}{2} )$.


To bound $\quot^{\prime\prime}(s)$, we recall from \eqref{quot-primeprime} that
\begin{align*}
	\quot^{\prime\prime}(s)
		&= i \sum_{k=0}^n
				\frac{ \nst_k e^{-is} }{(1 - \nst_k e^{-is})^2}
				\biggl( 1 - \frac{ 2 }{1 - \nst_k e^{-is}} \biggr)
\end{align*}
for $s \in ( -\frac{\delta}{2} , \frac{\delta}{2} )$.
Now, take a root $\nst \in \roots$ with polar representation $\nst = \nstrad e^{i \nstarg}$.
Then,
\begin{align}
	&(1 - \nst e^{-is})^2 e^{i(s-\nstarg)} \nonumber \\
	={}& e^{i(s-\nstarg)} - 2 \nstrad + \nstrad^2 e^{i(\nstarg-s)} \nonumber \\
	={}& (1+\nstrad^2) \cos(\nstarg-s) - 2 \nstrad - i (1-\nstrad^2) \sin(\nstarg-s) \nonumber\\
	={}&  -\frac12 (1+\nstrad)^2 \bigl( 1-\cos(\nstarg-s) \bigr)
		+\frac12 (1-\nstrad)^2 \bigl( 1+\cos(\nstarg-s) \bigr) \nonumber \\
	&	- i (1-\nstrad^2) \sin(\nstarg-s) \nonumber \\
	={}&  \frac12 \Bigl[
				(1+\nstrad) \sqrt{ 1-\cos(\nstarg-s) }
				- i \cdot \sgn_{\nstarg-s} \cdot
				(1-\nstrad) \sqrt{ 1+\cos(\nstarg-s) }
			\Bigr]^2,
		\label{quot-primeprime-exact}
\end{align}
where $\sgn_x \in \{ -1, 0, +1 \}$ denotes for $x \in \R$ the sign of $\sin(x)$.
Therefore,
\begin{align*}
	\vert 1 - \nst e^{-is} \vert
	&= \sqrt{ \bigl\vert (1 - \nst e^{-is})^2 e^{i(s-\nstarg)} \bigr\vert }
		\\
	&= \sqrt{
			\frac12 (1+\nstrad)^2 \bigl( 1-\cos(\nstarg-s) \bigr)
			+
			\frac12 (1-\nstrad)^2 \bigl( 1+\cos(\nstarg-s) \bigr)
		}
		\\
	&\geq \frac1{\sqrt{2}} (1+\nstrad) \sqrt{ 1-\cos(\nstarg-s) }
		\\
	&\geq \frac1{\sqrt{2}} (1+\nstrad) \frac{\sqrt{2}}{\pi} \bigl\vert \vert\nstarg\vert - \vert s \vert \bigr\vert
		\\
	&\geq \frac1{2\pi} (1+\nstrad) \vert\nstarg\vert,
\end{align*}
where we used that
$1-\cos(x) \geq \frac2{\pi^2} x^2$ for $x \in [-\pi, \pi]$ and that
$\bigl\vert \vert\nstarg\vert - \vert s \vert \bigr\vert \geq \frac{\vert\nstarg\vert}{2}$
since $\vert s \vert < \frac{\delta}{2} \leq \frac{\vert\nstarg\vert}{2}$.
Recalling the representation \eqref{quot-primeprime} of $q^{\prime\prime}$ we can now estimate separately both factors appearing there.
The first factor satisfies
\begin{align}
	\biggl\vert \frac{ -\nst e^{-is} }{(1 - \nst e^{-is})^2} \biggr\vert
	&\leq
	\frac{4\pi^2 \nstrad}{(1+\nstrad)^2\,\nstarg^2}
	\leq
	\frac{2\pi^2 \nstrad}{(1+\nstrad)^2\,\bigl( 1 - \cos(\nstarg) \bigr)}
	,
	\label{quot-primeprime-estimate}
\end{align}
where we used that $1-\cos(x) \leq \frac12 x^2$ for $x \in [-\pi, \pi]$.
For the second factor we get
\begin{align*}
	\biggl\vert 1 - \frac{ 2 }{1 - \nst e^{-is}} \biggr\vert
	&\leq
	1 + \frac{4\pi}{(1+\nstrad)\,\vert\nstarg\vert}
	\leq
	1 + \frac{4\pi}{(1+\nstrad)\,\delta}
	,
\end{align*}
where we used that $\delta \leq \vert\nstarg\vert$ by definition of $\delta$.
As an intermediate result, we get
\begin{align*}
	\vert \quot^{\prime\prime}(s) \vert
		&= \Bigg\vert
			i \sum_{k=0}^n
				\frac{ \nst_k e^{-is} }{(1 - \nst_k e^{-is})^2}
				\biggl( 1 - \frac{ 2 }{1 - \nst_k e^{-is}} \biggr)
			\Bigg\vert
		\\
		&\leq \sum_{k=0}^n
				\frac{2\pi^2 \nstrad_k}{(1+\nstrad_k)^2\,\bigl( 1 - \cos(\nstarg_k) \bigr)}
				\biggl( 1 + \frac{4\pi}{(1+\nstrad_k)\,\delta} \biggr)
\end{align*}
for $s \in \bigl( -\frac{\delta}{2} , \frac{\delta}{2} \bigr)$.
We emphasize that up to here all estimates hold for every possible root $\nst$, regardless of its location.
To proceed, we shall now specialize to the set-up of \ref{case:realrooted} or \ref{case:cyclotomic}.

\subsection{Specialization}

If a root $\nst$ of $\gf$ is real and negative or zero or if it satisfies $|\nst|=1$, then
\begin{align*}
	\biggl\vert \frac{-\nst}{(1-\nst)^2} \biggr\vert
	&=
	\frac{-\nst}{(1-\nst)^2},
\end{align*}
compare with \eqref{sigma-re}, \eqref{sigma-im} and with Figure~\ref{fig:sophoide}.
Further, in this case, \eqref{quot-primeprime-exact} for $s=0$ reduces to
\begin{align*}
	(1 - \nst)^2 e^{-i\nstarg}
	&=  \frac12 (1+\nstrad)^2 \bigl( 1-\cos(\nstarg) \bigr).
\end{align*}
Together with \eqref{quot-primeprime-estimate} this implies
\begin{align*}
	\biggl\vert \frac{ -\nst e^{-is} }{(1 - \nst e^{-is})^2} \biggr\vert
	&\leq
	\frac{2\pi^2 \nstrad}{(1+\nstrad)^2\,\bigl( 1 - \cos(\nstarg) \bigr)}
	=
	\pi^2 
	\biggl\vert \frac{-\nst}{(1-\nst)^2} \biggr\vert
	=
	\pi^2 \cdot
	\frac{-\nst}{(1-\nst)^2}
\end{align*}
so that
\begin{align*}
	\sup_{s \in ( -\frac{\delta}{2} , \frac{\delta}{2} )}
		\vert \quot^{\prime\prime}(s) \vert
	&\leq
	\sup_{s \in ( -\frac{\delta}{2} , \frac{\delta}{2} )}
		\Biggl\vert
			i \sum_{k=0}^n
				\frac{ \nst_k e^{-is} }{(1 - \nst_k e^{-is})^2}
				\biggl( 1 - \frac{ 2 }{1 - \nst_k e^{-is}} \biggr)
		\Biggr\vert
	\\
	&\leq
		\pi^2\sum_{k=0}^n \frac{-\nst_k}{(1-\nst_k)^2}
		\biggl( 1 + \frac{4\pi}{(1+\nstrad_k)\cdot\delta} \biggr)
	\\
	&\leq
		\sigma^2
		\, \frac{5\pi^3}{\delta}
	.
\end{align*}
Combining this with \eqref{alpha-estimate} yields
\begin{align*}
	\vert \alpha(t) \vert
	&\leq \frac{t^2}{2\sigma^3}
		\sup_{s \in ( -\frac{\delta}{2} , \frac{\delta}{2} )}
		\vert \quot^{\prime\prime}(s) \vert
	\leq \frac{5\pi^3}{2\sigma\delta} \, t^2
\end{align*}
for all $t \in ( -\frac{\sigma\delta}{2} , \frac{\sigma\delta}{2} )$.
Finally, \eqref{Stein-Tikho} leads to the bound
\begin{align}
	\sup_{w \in\R} \vert \PP(\widehat X \leq w) - F(w) \vert
	&\leq
		\frac{4}{3\sqrt\pi} \, \frac{5\pi^3}{2\sigma\delta}
		+ \frac{24}{\pi\sqrt{2\pi}}
		 \max\biggl\{ \frac{5\pi^3}{\sigma\delta} , \frac{2}{\sigma\delta} \biggr\}
		\nonumber\\
	&= \frac{c_0}{\sigma\delta},
		\label{main-result}
\end{align}
with a constant $c_0$ satisfying $c_0 := \frac{20\pi^3}{6\sqrt\pi} + 60\pi\sqrt{2\pi} < 531$. We are now in the position to complete the proof of Theorem~\ref{thm}.

\subsection{Proof of Theorem~\ref{thm} \ref{thm:i}}

First, we assume the set-up \ref{case:realrooted} of part \ref{thm:i} of Theorem~\ref{thm}.
Then, since all roots are real, we have that $\delta = \pi$.
The result thus follows from \eqref{main-result} with the constant $c_1$ given by $c_1:= \frac{c_0}{\pi} = \frac{20\pi^2}{6\sqrt\pi} + 60\sqrt{2\pi} < 169$.\qed

\subsection{Proof of Theorem~\ref{thm}~\ref{thm:ii}}\label{sec:proofii}

Next, we assume the set-up \ref{case:cyclotomic} of part \ref{thm:ii} of Theorem~\ref{thm}.
In this case, all roots $\nst\in\roots$ satisfy $|\nst|=1$ and $\delta$ is in relation with the fourth cumulant $\kappa$ of $X$.
Namely, in this case we get from \eqref{fourth-cumulant} that
\begin{align*}
	\kappa
	&= \sum_{k=0}^n
		\frac{ -\nst_k (1 + \nst_k )^2 - 2 \nst_k^2 }{(1 - \nst_k )^4}
		\\
	&= -\sum_{k=0}^n
		\frac{
			\bigl( \nst_k (1 + \nst_k )^2 + 2 \nst_k^2 \bigr)
			(1-\bar\nst_k)^4
		}{\vert 1 - \nst_k \vert^8}
		\\
	&= -\sum_{k=0}^n
		\frac{
			\bigl( \nst_k + 4 \nst_k^2 + \nst_k^3 \bigr)
						\bigl( 1 - 4\bar\nst_k + 6\bar\nst_k^2 - 4\bar\nst_k^3 + \bar\nst_k^4 \bigr)
		}{\vert 1 - \nst_k \vert^8}
		\\
	&= -\sum_{k=0}^n
		\frac{
			(\nst_k^3 + \bar\nst_k^3)
			- 9 (\nst_k + \bar\nst_k)
			+ 16
		}{\vert 1 - \nst_k \vert^8}
		\\
	&= -\sum_{k=0}^n
		\frac{
			2\cos(3\nstarg_k)
			- 18\cos(\nstarg_k)
			+ 16
		}{\vert 1 - \nst_k \vert^8}
	.
\end{align*}
From \eqref{quot-primeprime-exact} we know that
$\vert 1 - \nst \vert = \sqrt{2} \sqrt{1 - \cos(\nstarg)}$
for every root $\nst = e^{i\nstarg} \in \roots$.
Further, we have that $\cos(3x) = 4\cos(x)^3 - 3\cos(x)$ for all $x \in \R$.
This yields
\begin{align*}
	\kappa
	&= -\sum_{k=0}^n
		\frac{
			8\cos(\nstarg_k)^3
			- 24\cos(\nstarg_k)
			+ 16
		}{2^4 (1-\cos(\nstarg_k))^4}
		\\
	&= -\sum_{k=0}^n
		\frac{
			8(2+\cos(\nstarg_k))(1-\cos(\nstarg_k))^2
		}{16 (1-\cos(\nstarg_k))^4}
		\\
	&= -\sum_{k=0}^n
		\frac{
			2+\cos(\nstarg_k)
		}{2 (1-\cos(\nstarg_k))^2}
	.
\end{align*}
Note that if $n$ is even, this result already appeared as Equation (8) in \cite{HZ15}.
Finally, we use once again the fact that $1-\cos(x) \leq \frac12 x^2$ for $x \in [-\pi, \pi]$,
which implies the bound
\begin{align*}
	\vert\kappa\vert
	= -\kappa
	&\geq \sum_{k=0}^n
		\frac{2}{\nstarg_k^4}
	\geq \frac{2}{\delta^4}
	.
\end{align*}
Regarding the estimate \eqref{main-result}, this yields
\begin{align*}
	\sup_{w \in\R} \vert \PP(\widehat X \leq w) - F(w) \vert
	&\leq
		\frac{c_0}{\sigma\delta}
	\leq
		c_2\, \frac{\vert\kappa\vert^{1/4}}{\sigma}
	,
\end{align*}
with the constant $c_2$ given by  $c_2 := \frac{c_0}{2^{\frac14}} = \frac{10\pi^3 2^{\frac34} }{6\sqrt\pi} + 60 \pi^{\frac32} 2^{\frac14} < 447$. This completes the proof of Theorem \ref{thm}.\qed

\section*{Acknowledgment}

BR has been supported by the German Research Foundation (DFG) under project number 459731056. CT was supported by the DFG through SPP 2265 \textit{Random Geometric Systems} and CRC/TRR 191 \textit{Symplectic Structures in Algebra, Geometry and Dynamics}.


\begin{thebibliography}{999}

\bibitem{Bender}
E.A.~Bender:
Central and local limit theorems applied to asymptotic enumeration.
\emph{Journal of Combinatorial Theory. Series A} \textbf{15} (1973),
91--111.

\bibitem{BilleySwanson}
S.C.~Billey, J.P.~Swanson:
Cyclotomic generating functions.
\emph{arXiv: 2305.07620}.

\bibitem{Braenden}
P.~Br\"and\'en: Unimodality, log-concavity, real–rootedness and beyond.
In \emph{Handbook of Enumerative Combinatorics}. Chapman and Hall/CRC (2015).

\bibitem{CGS11}
L.H.Y.~Chen, L.~Goldstein, Q.-M.~Shao:
\emph{Normal Approximation by Stein's Method}.
Berlin and Heidelberg: Springer (2011).

\bibitem{Feller2}
W.~Feller:
\emph{An Introduction to Probability and its Applications, Vol.\ 2}.
Wiley (1966).

\bibitem{Harper}
L.H.~Harper:
Stirling behavior is asymptotically normal.
\emph{Annals of Mathematical Statistics} \textbf{38} (1967),
410--414.

\bibitem{HeertenEtAl}
N.~Heerten, H.~Sambale, C.~Th\"ale:
Probabilistic limit theorems induced by the zeros of polynomials.
To appear in \emph{Mathematische Nachrichten} (2023+).

\bibitem{HeilmannLieb}
O.J.~Heilmann, E.H.~Lieb:
Theory of monomer-dimer systems.
\textit{Communications in Mathematical Physics}, \textbf{25} (1972), 190--232.

\bibitem{HZ15}
H.-K.~Hwang, V.~Zacharovas:
Limit Distribution of the Coefficients of Polynomials With Only Unit Roots.
\emph{Random Structures \& Algorithms} \textbf{46} (2015),
707--738.

\bibitem{Lo61}
E.H.~Lockwood:
\emph{Book of Curves}.
Cambridge University Press (1961).

\bibitem{MichelenSara19Prep}
M.~Michelen, J.~Sahasrabudhe:
Central limit theorems and the geometry of polynomials.
\emph{arXiv: 1908.09020}.

\bibitem{Pitman}
J.~Pitman: 
Probabilistic bounds on the coefficients of polynomials with only real zeros.
\emph{Journal of Combinatorial Theory. Series A} \textbf{77} (1997),
279--303.

\bibitem{Ro22}
A.~R\"ollin:
Kolmogorov bounds for the normal approximation of the number of triangles in the Erd\H{o}s-R\'enyi random graph.
\emph{Probability in the Engineering and Informational Sciences} \textbf{36} (2022),
747--773.

\bibitem{Ti80}
A.N.~Tikhomirov:
On the rate of convergence in the central limit theorem for weakly dependent random variables.
\emph{Theory of Probability and Its Applications} \textbf{25} (1980),
790--809.

\end{thebibliography}
\end{document}